
\magnification=\magstephalf \documentstyle {amsppt}

\def\a{\alpha}

\def\b{\beta}

\def\g{\gamma}

\def\G{\Gamma}

\def\d{\delta}

\def\D{\Delta}

\def\e{\epsilon}

\def\z{\zeta}

\def\k{\kappa}

\def\l{\lambda}

\def\L{\Lambda}

\def\m{\mu}

\def\n{\nu}

\def\c{\chi}

\def\si{\sigma}

\def\Si{\Sigma}

\def\o{\omega}

\def\bz{\Bbb Z}
\def\bc{\Bbb C}
\def\({\left(}
\def\){\right)}
\def\]{\right]}
\def\[{\left[}

\hfuzz 100pt \topmatter \title Skein theory and Witten-Reshetikhin-Turaev
Invariants of links in lens spaces \endtitle \rightheadtext{Skein theory and
WRT
Invariants}

\author Patrick M. Gilmer \endauthor \affil Louisiana State University
\endaffil
\address Department of Mathematics,  Baton Rouge, LA 70803 U.S.A  \endaddress
\email gilmer\@math.lsu.edu \endemail

\abstract We  study  the  behavior of the Witten-Reshetikhin-Turaev $SU(2)$
invariants  of an arbitrary link in $L(p,q)$ as a function of the level $r-2$.
They are given
by $\frac{1}{\sqrt{r}}$ times one of $p$ Laurent polynomials evaluated at
$e^{\frac {2 \pi
i} {4pr}}.$ The congruence class of $r$ modulo $p$ determines which polynomial
is
applicable.  If $p \equiv 0 \pmod{4},$ the meridian of $L(p,q)$ is non-trivial
in
the skein module but has trivial Witten-Reshetikhin-Turaev $SU(2)$ invariants.
On
the other hand, we show that one may recover the element in the Kauffman
bracket
skein module of $L(p,q)$ represented by a link from the collection of the WRT
invariants at all levels if $p$ is a prime or twice an odd prime. By a more
delicate argument, this is also shown to be true for $p=9.$  \endabstract

\thanks This research was partially supported by a grant from the Louisiana
Education Quality Support Fund \endthanks \keywords  skein,  lens space, Gauss
Sum \endkeywords \subjclass 57M99 \endsubjclass

\endtopmatter \document

\centerline{This version: 7\slash 6\slash 98, First Version:1\slash 28\slash
98} \head  \S 1 Introduction    \endhead

We consider the Witten-Reshetikhin-Turaev $SU(2)$ invariants
$\{w_r(M,J)\}_{r \ge   2} \in \Bbb C.$  \cite{W,RT} These are the invariants of
Witten at level $r-2.$  We specify the precise version of the WRT invariant
that we study near the  beginning of \S2.

Let $M$ be a closed oriented 3-manifold. Let $S(M)$ be the Kauffman bracket
skein
module. Recall this is a module over $\L= \Bbb Z[A,A^{-1}].$ Elements of $S(M)$
are
linear combinations of links over $\L$ modulo the Kauffman
relations. Let $J$ be a framed link in $M,$ then $J$ represents an element
in $S(M),$ called its skein class.  We will use the symbol $\Cal J$ to
represent a general skein element.

It is immediate from the definition of $w_r$ that $w_r(M,J)$ only depends on
the skein class of J.
In fact we may extend the definition of $w_r$ to  $S(M),$ by extending
linearly, and sending $A$ to
$-a_r.$  Here and throughout, we adopt the notation that $\xi_g=e
^{\frac {2 \pi i}{ g}},$ and $a_r=\xi_{4r}.$
We may study the behavior of $w_r$ on links in $M$ by
calculating  $w_r$ for generators of $S(M).$

Hoste and Przytycki \cite{HP} have calculated the skein module for $L(p,q).$
Throughout this paper $p$ will be an integer greater than one. $S(L(p,q))$ is
the
free  $\L$ module generated by $[\frac p 2]+1$ generators: $x_0,
 x_1, x_2 \cdots   x_{[\frac p 2]}.$ Here $x_0$ is given by the empty link.
We  denote this by module by $\L(x_0, x_1 \cdots  x_{[\frac p 2]}).$ When
$L(p,q)$ is presented as $-p\slash q$ surgery to the unknot, $x_1$ is given by
a
meridian to this unknot, and  $x_i$ is given by $i$ parallel  meridians. We may
consider $\m_c,$ a meridian colored by $c,$ as an element in
$S(L(p,q)).$ By this, we mean one should place the skein element $e_{c} \in
S(S^1
\times D^2)$  defined in \cite{BHMV1} into a tubular neighborhood of a
meridian.
The recursion definition in $S(S^1 \times D^2)$ for $e_c$  is $e_c= \a e_{c-1}-
e_{c-2}, \quad
e_0=1, e_{-1}=0.$   It is useful to define
$e_c$ for $c$ negative as well,  by running the above recursion formula
backwards.
The recursion definition implies that  $\a^c \in \L(e_0,\cdots ,e_c),$ for
$c\ge 0.$
Here we let $\a$ denote $S^1 \times {0}$ with the standard framing in $S^1
\times D^2.$
Thus
$S(L(p,q))$  is also the $\L$ module generated by : $\m_0,   \m_1,
\m_2 \cdots   \m_{[\frac p 2]}.$

We wish to study $w_r( L(p,q),\m_{c})$ for $r\ge 2,$ and $p \ge 2.$  We will be
particularly interested in $0 \le c \le  [\frac p 2].$

The answer involves generalized Gauss sums:

$$\Cal G_p(a,b)= \sum_{n=0}^{p-1}  \xi_p^{a n^2+bn}.$$

Section two of \cite{LL} gives two convenient theorems which allow the
evaluation of such sums. We will refer these theorems frequently.
For $k
\in \bz,$ consider the following generalized Gauss sum: $$G_\pm(p,q,c,k)=
\Cal G_p(qk ,qc+q \pm1).$$

Note that $$G_\pm(p,q,c,k)=G_\pm(p,q,c,k+p)=G_\pm(p,q,c+p,k)=
 \overline{G_\pm(p,q,c,-k)}.$$

For $ 0 \le k < p,$ define  $f_{p,q,c,k}(z)$ to be the Laurent polynomial in
$\bc [z,z^{-1}]$ :
$$\frac {i (-1)^{c+1}}{\sqrt{2p}} z^{12 p\ s(q,p)  +q (c^2+2c)}
\left( G_+(p,q,c,k)\ z^{2(c+1)}
- G_-(p,q,c,k)\ z^{-2(c+1)} \right) .$$

Here $s(q,p)$ is a Dedekind sum. Then
$$f_{p,q,c,k}(z)=- \overline{f_{p,q,c,p-k}(\bar z)}.$$

The following theorem in the case $c$ is zero is Jeffrey's formula \cite{J} for
$w_r( L(p,q)).$ In \S2, we will prove this theorem by adapting Jeffrey's
proof.

\proclaim{Theorem 1} If $r \equiv  k \pmod{p},$ $$\sqrt{r}\
w_r( L(p,q),\m_{c})= f_{p,q,c,k}(\xi_{4pr}) $$ \endproclaim

Suppose $J$ is a link in $L(p,q)$, let $ [J]=\sum_{c=0}^{[\frac p 2]}
C_{J,c}(A)
[\m_{c}] \in S(L(p,q)).$ Note $C_{J,c}(A) \in \L.$ Then  $\sqrt{r}\
w_r( L(p,q),J)= \sum_{c=0}^{[\frac p 2]} f_{p,q,c,r}(\xi_{4pr})
C_{J,c}(-a_r).$ For $ 0 \le k <p,$ let $$f_{J,k}(z)= \sum_{c=0}^{[\frac p 2]}
f_{p,q,c,k}(z) C_{J,c}(-z^{p}) \in \bc[z,z^{-1}] \leqno{(1.1)}$$

Since $C_{J,c}$ has integral and therefore real coefficients,
$$f_{J,k}(z)=- \overline{f_{J,p-k}(\bar z)}\leqno{(1.2)}$$

\proclaim {Corollary 1 } If $J$ is a link in $L(p,q)$ (or any element of
$S(L(p,q))$),
then for each  $ 0 \le k <p,$ there is a  polynomial $f_{J,k}(z)$ in
$ \bc[z,z^{-1}]$ such that
for $r \equiv k \pmod{p},$ $$\sqrt{r}\ w_r( L(p,q),J)= f_{J,k} (
\xi_{4pr}).$$  \endproclaim

Thus there are $p$ Laurent  polynomials which describe $\sqrt{r} \
w_r( L(p,q),J)$ as a function of $ \xi_{4pr}.$ Each  function is applicable
if one
restricts  $r$ to  a particular congruence class modulo $p.$ Moreover the last
$p-[\frac p 2 ]-1$ of these polynomials are determined in a simple way by the
first $[\frac p 2 ]+1$ of these polynomials.

In \cite{G1}, we showed that if $J$ is a link in a connected sum of $g$  $S^1
\times S^2$ 's, then   $w_r( S^3)^{g-1}\  w_r(\sharp^g S^1 \times S^2,
J)$ is given by a rational function in $a_r$ for $r$  large enough (depending
on $J$).  We also gave examples where this formula fails for small $r.$
Moreover if $g=0$ or one, this rational function is a Laurent  polynomial.
We were interested in whether similar statements could be made about links in
lens spaces. This was the genesis of the present paper.

Lawrence and Rozansky have described the Witten-Reshetikhin-Turaev $SU(2)$
invariants of a Seifert rational homology sphere by a single holomorphic
function of $r$ \cite{LR}. See also \cite{L1,L2}. They also show that Otsuki's
power series in $h=q-1$ evaluated by setting q to be an rth root of unity
converges r-adically to the  $r$-th Witten-Reshetikhin-Turaev invariant
(associated to $SO(3)$), if $r$ is an odd prime which is relatively prime to
the order of the first homology.  Suppose we normalize  $w_r$  to be 1 on $S^3$
by dividing by $w_r( S^3)$. Then $w_r( S^3)^{-1}w_r(  L(p,q),J)$
for fixed congruence class of $r$ modulo $p$ is given by a rational function of
$ \xi_{4pr}.$  As we vary the congruence class of $r$ modulo $p,$ we obtain
different functions. Note the case $c=0$ in Theorem 1, which is really just a
reinterpretation of
Jeffrey's result \cite{J, \rm Theorem (3.4)}, concerns $L(p,q)$ with the empty
link. Thus work of Lawrence and Rozansky  applies in this case, but describes
the dependence on $r$ in a different way. We use several rational functions of
$ \xi_{4pr}.$ Lawrence and Rozansky use a single  holomorphic function of $r.$
This holomorphic function is defined by a line integral.  Recently  the r-adic
convergence mentioned above has been extended to links in a general rational
homology sphere \cite{R}.

Using \cite{LL}, one may see that if $p \equiv 0 \pmod{4},$ and $c$ is odd,
then
$G_\pm(p,q,c,k)$ is zero.

\proclaim{Corollary 2} If $p \equiv 0 \pmod{4}$ and $c$ is odd,
$w_r( L(p,q),\m_c)=0.$ Similarly, if $p \equiv 0 \pmod{4}$ and $c$ is
odd and positive, $w_r( L(p,q),x_c)=0.$ \endproclaim

To see the last statement, we note: If $c$ is odd and positive,  $\a^c \in
\L(e_1,e_3\cdots
,e_c)\subset S(S^1 \times D^2).$  In particular for any $n \ge 1$, the meridian
in $L(4n,q)$  represents a
non-trivial skein class but has all WRT invariants trivial.

On the other hand, we can sometimes show that the skein class of a link may be
recovered from its WRT invariants. Suppose we are given $\sqrt{r}\ w_r(
L(p,q),J)$ for all $r,$ then $f_{J,k}$ is determined for all $k$ in the
range $ 0\le k <p .$ This is because a non-zero Laurent polynomial cannot have
infinitely many roots.
Let $\Cal R$ denote the field of rational functions in $z$ with complex
coefficients. Suppose the  $p \times \([\frac p 2]+1\)$ matrix
$[f_{p,q,c,k}(z)]_{0\le k \le p,0\le c \le [\frac p 2]}$ over $\Cal R$
has rank $1+[\frac p 2].$ Then we could solve the set of equations (1.1) for
$C_{J,c}(-z^p).$ Thus
$C_{J,c}(A)$ would be determined for all $c.$

Let $\Bbb S(L(p,q))$ denote the skein module over $\Cal R$ where $A= -z^{p}.$
This is the free $\Cal R$ module generated by framed links in $\Bbb S(L(p,q))$
modulo the submodule generated by isotopy and the usual Kauffman relations
where now we take  $A$ to be  $-z^{p}.$ We may then  define $w_r$ for skein
classes in $\Bbb S(L(p,q))$ by sending z to $\xi_{4pr}.$ Equation (1.1) will
still hold but Equation (1.2) will not, as $C_{J,c}(z)$ need not have real
coefficients.
If   $[f_{p,q,c,k}(z)]_{0\le k \le p,0\le c \le [\frac p 2]}$ over $\Cal R$ has
rank less than $1+[\frac p 2],$ then there will be non-trivial classes in $\Bbb
S(L(p,q))$
with trivial WRT invariants.
In \S 3 and \S4, we will study the rank of this matrix to prove:

\proclaim{Theorem 2}   Skein classes
$\Cal J \in \Bbb S (L(p,q))$ are determined by  $w_r( L(p,q),\Cal J)$ for all
$r$ if and only if $p$ is a prime or twice an odd prime \endproclaim

In particular, there is a non-trivial skein class $\Cal J$ in $\Bbb
S(L(9,q))$ with trivial WRT-invariants. The subset of skein classes with
trivial
WRT-invariants turn out to form a one dimensional line in the $\Cal R$-vector
space $\Bbb S(L(9,q)).$ Moreover this line does not intersect the ordinary
skein
module $S(L(9,q)),$ except at the trivial class. In this way, we will prove in
 \S5:

\proclaim{Theorem 3} Let  $q$ be prime to 9. Every skein class $\Cal J$ in
$S(L(9,q))$
is determined by  $w_r(  L(9,q),\Cal J)$ for all $r.$\endproclaim

\head \S2 Proof of Theorem 1 \endhead

This result is obtained by studying Jeffrey's calculation of
the WRT invariant of lens spaces.  However we must first relate this
calculation to the
skein theoretic definition of the invariant.

We  find it convenient to use  the
variant of the \cite{BHMV2} TQFT  defined in \cite{G2}, where $p_1$ structures
are
replaced, following Walker \cite{Wa}, with Lagrangian subspaces and integers.
A
closed oriented
surface  (possibly empty) $\Si$  together with a  choice of Lagrangian subspace
 of $ H_1(\Si,\Bbb Z)$  is an object of a cobordism category    $\Cal C_{2r}.$
A  morphism from $\Si$ to $\Si'$ is a compact oriented
3-manifold $M,$ weighted by an integer $w,$ and containing a (possibly empty)
colored framed link $J$ such that  boundary of $M$ comes equipped with an
identification with $-\Si
\coprod \Si'.$   For simplicity, we have only describe some of the objects and
morphism in $\Cal C_{2r}.$
The TQFT assigns modules $ V_{2r} (\Si)$ to objects $\Si.$ A morphism $(M,J,w)$
from $\Si$ to $\Si',$ defines a linear map $Z_{2r}(M,J,w): V_{2r} (\Si)
\rightarrow  V_{2r} (\Si)'.$ $V_{2r}(\emptyset)$ is the ground ring $k_{2r}.$
If $M$ has no boundary, $Z_{2r}(M,J,w)$ is multiplication by $<(M,J,w)>_{2r}.$

We can now specify the version of the WRT SU(2) invariant that we study.
$w_r(M,J)= \psi (<(M,J,0)>_{2r}).$  Here $\psi:k_{2r}
\rightarrow \bc$ is the homomorphism \cite{MR, note p.134} which sends $A$ to
$-a_r,$ and sends $\k$ to $\z a_r^{-1}.$ Here   $\z= \xi_8.$

$L(p,q)$ is obtained by gluing together two copies of $S^1 \times D^2$ by the
self-homeomorphism $U$ of $S^1 \times S^1$ given with respect to a
meridian-longitude basis by the matrix $(\smallmatrix q & b \\ p&
d\endsmallmatrix)\in SL(2,\bz).$ $L(p,q)$ weighted by an integer $n$ as
a morphism
of the cobordism category is denoted $(L(p,q),n).$ Let $\Si$
denote the torus $S^1 \times S^1$ assigned the Lagrangian subspace generated by
the meridian $[S^1 \times {1}].$  For each integer $c,$ let $H_c$
be the morphism from the empty set to $\Si$ given by $S^1 \times D^2 $ weighted
zero, with the core $S^1 \times {0} $ with the standard framing colored  by
$c.$ Let $\tilde H_c$ be this same manifold (with the opposite orientation)
viewed as a morphism from  $\Si$ to the empty set.
In $\Cal C_{2r},$ we have that $(L(p,q),0)$ is the composition $ \tilde H_0
\circ C(U,0)
\circ H_0.$ We let $C(U,w)$ denotes the mapping cylinder of $U$ weighted by an
integer $w.$  More
generally, we  have that  $L(p,q)$ weighted zero with the meridian colored $c$
is
given by the composition  $\tilde H_c \circ C(U,0) \circ H_0.$ $H_c$ represents
an
element $\b_c$ in $V_{2r}(\Si).$ As $c$ ranges over $ 0 \le c \le r-2,$ these
elements  form an orthonormal basis $\Cal
B.$ We will also need elements $\n_s= (-1)^{s-1}\b_{s-1}.$ By \cite{BHMV1,6.3},
they have the following symmetries
$$\n_l= \n_{2r+l} \quad \text{and} \quad
\n_l= -\n_{-l}\quad\forall l\in \Bbb Z \leqno{2.1}$$
Note that $w_r( L(p,q),(-1)^{l-1} \mu_{l-1})$ will possess the same symmetries.
Let  $\Bbb B$ denote the basis with elements $\n_s$:  for $1\le s\le r-1.$

As usual,  we let $S=(\smallmatrix
0&-1\\1 &0\endsmallmatrix),$  and $T=(\smallmatrix 1&1\\0 &1\endsmallmatrix).$
Let $J(m)= T^m S.$ Then $C(J(m),0)=C(T,0)^m \circ C(S,0),$ as the contribution
of
$\si$ (Maslov index) in the formula for the computation of the weight of  a
composition is zero. Following \cite{J} we expand $U= J(m_t) \circ J(m_{t-1})
\cdots J(m_1),$ where $m_t=0,$ and $t>1.$ Let $U_i= J(m_i) \cdots J(m_1).$ Let
$c_i$ denote
the lower left entry in $U_i.$  For $1\le i\le t,$ define $w_i$ by $ C(U_i,w_i)
=C(J(m_i),0)\circ \cdots  \circ C(J(m_1),0).$ So $w_1=0.$ For $1< i\le t,$ we
have $ C(U_i,w_i)
=C(J(m_i),0)\circ C(U_{i-1},w_{i-1}).$ Computing the contribution of $\si,$ one
has that $w_i= w_{i-1}+\text{sign}(c_{i-1}c_i).$ Thus $w_t =  \sum_{i=2}^{t}
\text{sign}(c_{i-1}c_i).$  This is the signature of the matrix $W_L$ of
\cite{J,Equation 3.6}, at least assuming this matrix is a definite matrix which
we may
insist  without loss of generality. Note that
$\text{Trace} (W_L)=\sum_i m_i.$

With respect to $\Bbb B,$ $Z_{2r}(T,0)$ is given by $ \hat T_{l,j}=\d_l^j
(-A)^{l^2-1},$ and $Z_{2r}(S,0)$ is given by $\hat S_{l,j} =\eta_{2r} [l j] $
\cite{G2}. Here $j$ and $l$ range from $1$ to $r-1.$ Also $[k]$ denotes the
quantum integer $\frac {A^{2k} - A^{-2k}} {A^{2} - A^{-2}},$ and
  $\eta_{2r}$ is the scalar of \cite{BHMV2}  with $\psi(\eta_{2r})= - \frac
{a_r^2-a_r^{-2}}{\sqrt{2r}} i.$ Witten \cite{W} used a unitary matrix
representation $\Cal
R_r$ of $SL(2,\bz),$ with $\Cal R_r
(S)=\psi (\hat S)$ and
$ \Cal R_r (T)= a_r \z^{-1}\psi (\hat T) = \d_l^j \z a_r^{l^2}.$
In \cite{J}, Jeffrey found an explicit formula for $\Cal
R_r \left([\smallmatrix a&b\\c& d \endsmallmatrix]\right)$ in terms of $a,$
$b,$ $c,$ and $d.$
With respect to $\Bbb B,$ we have
$$(\z a_r^{-1})^{3
\text{Sign}(W_L)}\psi\left(  Z_{2r}(C(U,0))\right)= \psi (Z_{2r}(C(U,w_t)))=
(\z
a_r^{-1})^{\text{Trace}(W_L)} \Cal R_r(U).$$
Since  the Rademacher $\Phi$-function
is given by \cite{J,3.2}\cite{KM} $$\Phi (U)=\text{Trace}(W_L)-3 \
\text{Sign}(W_L),$$ we have $$ \psi \left( Z_{2r}(C(U,0))\right) =  \left( \z
a_r^{-1}\right) ^{\Phi
(U)} \Cal R_r(U).$$

This is exactly the  correction factor Jeffrey used to find the Witten
invariant
of $L(p,q)$ in the canonical framing\cite{J, Lemma 3.3, Theorem 3.4}. We have
for $1\le l \le r-1:$
$$w_r( L(p,q),(-1)^{l-1} \mu_{l-1})=
\psi\left(<\nu_l,Z_r (C(U,0)) \nu_1 >\right)=  (\z
a_r^{-1})^{\Phi (U)} \Cal R_r(U)_{l,1}$$

Jeffrey \cite{J,Theorem 3.4} derived and simplified  an expression for  $(\z
a_r^{-1})^{\Phi (U)} \Cal R_r(U)_{1,1}.$ One may adapt her calculation as
follows. Extending \cite{J, Equation (3.8$\frac 1 2$)} to the case $l \ne 1,$
we have in our notation

$$w_r( L(p,q),(-1)^{l-1} \mu_{l-1})
 = \frac { -i
}{\sqrt{2rp}} a_r^{-\Phi(U)}\xi_{4rq}^b \sum_{n=1}^{p} \sum_{\pm}\pm
\xi_{4rpq}^{(q \g \pm 1)^2} \leqno {2.2}$$ where $\g=l+2 rn,$ and
$1 \le l\le r-1.$

Note  that $q \g^2 \pm 2\g$ modulo $4rp$ only depends on $n$  modulo $p.$  Let
$S_\pm(l)$ denote the unordered list, with multiplicity, of the values of
$q \g^2 \pm 2\g$ modulo $4rp$ as $n$ varies from 1 to $p.$
Then $S_\pm(l+2r)=S_\pm(l),$ and $S_\pm(-l)=S_\mp(l).$ So $\sum_{n=1}^{p}
\sum_{\pm}\pm
\xi_{4rpq}^{(q \g \pm 1)^2}= \xi_{4rpq} \sum_{n=1}^{p} \sum_{\pm}\pm
\xi_{4rp}^{q \g^2 \pm 2\g}$ remains the same when $l$ is replaced by $l+2r,$
and changes sign when $l$ is replaced by $ -l.$
By (2.1), it follows that (2.2) holds for all integers $l.$

As $$(q \g \pm
1)^2= q^2l^2+1 \pm 2ql + 4rq(qrn^2+qln \pm n),$$ we have: $$ \align
w_r( L(p,q),\mu_c) = \frac { i (-1)^{c+1}} {\sqrt{2rp}} \sum_{\pm}\pm
\xi_{4rpq}^{pb - pq \Phi(U) + q^2 l^2 +1 \pm 2q l}
\sum_{n=1}^{p}\xi_p^{(qr)n^2+(q l \pm 1)n}\\ = \frac { i (-1)^{c+1}}
{\sqrt{2rp}}
\xi_{4rpq}^{pb - pq \Phi(U) + q^2 l^2 +1}\sum_{\pm}\pm  \xi_{2rp}^{\pm  l}
G_{\pm}(p,q,c,r). \endalign$$

Here $l$ denotes $c+1.$ This holds for all integers $c.$
But now making use of $1= \det(U)= qd-bp,$ and the definition of $\Phi,$ we
have:
$$pb - pq \Phi(U) + q^2 l^2 +1=q^2(l^2-1)+12 pq\ s(d,p)=q^2(l^2-1)+12 pq \
s(q,p)$$

\noindent As in \cite{J}, we have noted that $d\equiv q^*\pmod{p},$ where $q^*$
is an integer with $q q^* \equiv 1 \pmod{p},$ and that
$s(d,q)=s(q^*,p)=s(q,p).$
Thus

$$ \sqrt{r} \ w_r( L(p,q),\mu_c ) = \frac { i (-1)^{c+1}}{\sqrt{2p}} {}
\xi_{4rp}^{12 p\ s(q,p)  +q(c^2+2c)}\sum_{\pm}\pm \xi_{4rp}^{\pm 2(c + 1)}
G_{\pm}
(p,q,c,r) $$

A key observation for this paper is that $r$ enters the right hand side of this
last formula in only two places in $\xi_{4rp}$ which we may think of as a
variable $z$ and in $G_{\pm} (p,q,c,r)$ where its contribution only depends on
$r
\pmod{p}.$

Thus  if $r \equiv  k \pmod{p},$ then
$$\sqrt{r} \ w_r( L(p,q),\m_{c}) = f_{p,q,c,k}(\xi_{4pr}).$$

\head \S 3 Proof of theorem 2 (only if part)  \endhead

Using \cite{LL} one may prove:

\proclaim{Lemma 1} If $a \not\equiv 0 \pmod{p},$ and $b^2\equiv  {b'}^2
\pmod{p},$
then $\Cal G_p(a,b)=\Cal G_p(a,b').$ \endproclaim

Let $\#_p$ denote the number of squares modulo $p.$

\proclaim{Lemma 2} The
number of distinct columns appearing in either of the matrices
$[G_+(p,q,c,k)]_{1\le k < p, 0\le c \le [\frac p 2 ]}$
or $[G_-(p,q,c,k)]_{1\le k < p, 0\le c \le [\frac p 2 ]}$ is less than  or
equal to $\#_p.$ \endproclaim

\proclaim{Lemma 3} The rank of  $[f_{p,q,c,k}(z)]_{0\le k < p,0\le c \le [\frac
p 2]}$ is less than or equal to $1+\#_p.$\endproclaim

\demo{Proof} The columns of $[f_{p,q,c,k}(z)]_{1\le k < p,0\le c \le [\frac p
2]}$ are linear combinations of the columns discussed in Lemma 2. Thus the rank
of this matrix is less than or equal to $\#_p.$  $[f_{p,q,c,k}(z)]_{0\le k <
p,0\le c \le [\frac p 2]}$ has one extra row.  \qed \enddemo

\proclaim{Lemma 4} If $p$ is not prime, twice an odd prime, four or nine, then
$\#_p <[\frac p 2 ].$ \endproclaim

\demo{Proof}  Suppose $(a,b)=1,$  then
$\bz_{ab} = \bz_a \oplus \bz_b $ as rings, and  so $\#_{ab} =\#_{a}  \#_{b}.$
If $\#_{a} < [\frac a 2],$  then the lemma holds for $p= a b,$ where $(a,b)=1.$
 If $p$ satisfies the hypothesis, then it can be written as $a \cdot b$ with
$(a,b)=1,$ and $a$ of one of the following forms. In each case, we show
why  $\#_{a} < [\frac a 2].$

(i) $a=2^t$ with $t \ge 3.$ The case $a=8$ may be easily checked, so suppose
$t\ge 4.$  Then every square is the square of a number in the range from zero
to $2^{t-1}.$ However zero, $2^{t-2},$ and $2^{t-1}$ are in this range and have
the same square.

(ii) $a=s^{2t}$ with $s$ an odd prime and $t \ge 1.$ If $t=1$, assume $s \ge
5,$ i.e. $a \ne 9.$
Then every square is the square of a number in the range from zero to $[\frac a
2].$  However zero, $s^t,$ and $2s^t,$ are in this range and have the same
square.

(iii) $a=s^{2t+1}$ with $s$ an odd prime, and $t \ge 1.$
Then every square is the square of a number in the range from zero to $[\frac a
2].$  However $0^2= (s^{t+1})^2,$ and
$\left(\frac {s^{t+1} +s^{t}}{2}\right)^2=\left(\frac {s^{t+1}
-s^{t}}{2}\right)^2$ are the squares of numbers from this range.

(iv) $a= s_1 \cdot s_2,$ where $s_1$ and $s_2$
are distinct odd primes. So   $\#_{a}=\#_{s_1 \cdot s_2}=
\#_{s_1} \cdot \#_{s_2}=
\frac {(s_1+1)} 2
\cdot \frac {(s_2+1)} 2
< \left[\frac {s_1 \cdot s_2} 2\right]= \left[\frac a 2\right].$

(v) $a$ is four times an odd prime,  nine times an odd prime (other than 3) ,
$2 \cdot  9,$ or  $4 \cdot  9.$ These cases are dealt with as in (iv).
\qed \enddemo

These lemmas yield the only if part of Theorem 2, except for the cases $p=4,$
and $p=9.$ When $p=4,$ $[f_{p,q,c,k}(z)]_{0\le k < p,0\le c \le [\frac p 2]}$
has a column of zeros by the proof of Corollary 2. When $p=9$ direct
calculation (see \S 5) shows this matrix has rank four.

\head \S 4 Proof of theorem 2 (if part)  \endhead

Let $\hat c$ be defined by $0\le \hat c \le p-1$ and $\hat c \equiv
q_p^*(c-1)-1\pmod{p}.$  Note that  $q\hat c +q+1\equiv c \pmod{p}.$
So that:

$$G_+(p,q,\hat c,k)=  \Cal G_p(qk,c) $$

Let $\D$ and $\G$ be  functions from $\{ 0, \cdots [\frac p 2]\}$ to
$\{ 0,\cdots p-1\}$  such that $\D(0)=0,$ and  $\G(0)=\hat 0.$
We will show that $[f_{p,q,c,k}(z)]_{0\le k < p,0\le c \le [\frac p 2]}$
has rank $1+[\frac p 2]$ by showing the determinant
$[f_{p,q,\G(c),\D(k)}(z)]_{0\le k \le [\frac p 2],0\le c \le [\frac p 2]}$
 is non-zero for some choice of $\G$
and $\D.$
It suffices to show the determinant of $[G_+(p,q,\G(c),\D(k))(z)]_{0\le k \le
[\frac p 2],0\le c \le [\frac p 2]}$ is nonzero.

Note $$G_+(p,q,\hat c,0) = \Cal G_p(0,c)= \cases p,&   \text{if $c=0$}\\
0,&\text{ if $c \ne
0$}.\endcases $$
Thus we only need to show $[G_+(p,q,\Delta (c),\G(k))]_{1\le k \le [\frac p 2],
1\le c \le [\frac p 2]}$ is nonsingular over $\Bbb C,$ for some $\D$ and $\G.$
If $p=2,$ this matrix is a one by one matrix with entry $\Cal G_2(1,1)=2.$
Thus we may assume $p$ is an odd prime or twice an odd prime.

\subhead {\bf The case $p$ is  an odd prime} \endsubhead

Let $t_u^*$ be defined by $t t_u^*\equiv 1 \pmod{u},$ if $t$ is prime to $u.$
Define $\D(k)= k_p^*,$ if $ k\ne 0.$ Define $\G (c) = \hat c.$
We note that for $k \ne 0,$ we have by \cite{LL}

$$ G_+(p,q,\hat c, k^*_p )=\Cal G_p(q k_p^*,c)=\e (p) \left({\frac
{qk^*_p}{p}}\right)\ \sqrt{p}\  \left(_q\l_{p}^{c^2}\right)^k $$

Here $_q\l_{p}$ denotes the  primitive $p$th root of unity $\xi_p^{-q^*(\frac
{p+1}{2})^2}.$ Also $\e(p)$ is one or $i$ accordingly as $p$ is one or
three modulo four, and $\left({\frac {a}{b}}\right)$ denotes a Jacobi symbol.
So the rows of $\det[G_+(p,q,\hat c,k^*)]_{1\le k \le \frac p 2,
1\le c \le \frac p 2 }$  are non-zero  multiples of the rows  of a Vandermonde
matrix with a non-zero determinant. Here we
make use of the fact that the numbers $_q\l_{p}^{c^2},$ as $c$ varies from $1$
to
$[\frac {p}{2}]$  are all distinct. For this it is important that as $c$ ranges
over $1\le c \le [\frac p 2],$ $c^2$ does not repeat itself modulo $p.$

\subhead {\bf The case $p$ is twice an odd prime $s$} \endsubhead

For this case, we define $\G$ by saying that $\G (0),$ $\G (1) ,$ $\cdots,$ $\G
(s)$ in sequence are given by $\hat 0,$ $\hat 2,$ $\hat 4,$ $\cdots,$
$\hat{(s-1)},$$\hat 1,$$ \hat 3,$$  \cdots$$ \hat s.$  Also define  $\D$ by
saying $\D (0),$$ \D (1) ,$$\cdots,$ $\D (s)$ in sequence are given by $0,$
$2(1^*_s),$ $2(2^*_s),$ $\cdots,$ $2( \frac{s-1}{2} )_s^*,$ $ 1^*_p,$ $ 3^*_p,
$
$\cdots,$ $(s-2)^*_p,$ $s.$  Since $p$ is even, $q$ is odd, and $\hat c \equiv
c
\pmod{2}.$ Thus for both $\G$ and $\D,$  even values come first, and then  odd
values follow.

Next notice that by \cite{LL} $$ G_+(2s,q,\hat c ,s)=\Cal G_{2s}(q s,c)= \cases
0,&
  \text{if $c$ is relatively prime to $s$ }\\ 2s,&\text{ if $c =s$}.\endcases
$$

\noindent Thus the last row  of $[G_+(p,q,\Gamma (c),\Delta(k))]_{1\le k \le
s, 1\le c \le s}$ consists of all zeros except for the last entry of $2s.$ Thus
we only need
to show $[G_+(p,q,\Gamma (c),\Delta(k))]_{1\le k \le s-1, 1\le c \le s-1}$ has
a
nonzero determinant. For $j$ relatively prime to $s,$  and $c \not \equiv j
\pmod {2}$ by \cite{LL}: $$ G_+(2s,q,\hat c ,j)=\Cal G_{2s}(q j,c)= 0 $$
Thus our matrix is a block matrix of the form $[\smallmatrix A&0\\0&B
\endsmallmatrix ].$

We now study the matrix $A.$ For $1 \le d \le \frac{s-1}{2},$ and $1 \le j \le
\frac{s-1}{2},$ by \cite{LL} we have:

$$ G_+(2s,q,\hat {(2 d)}  ,2 j^*_s)=\Cal G_{2s}(2q j^*_s,2d)= 2 \Cal G_{s}(q
j^*_s,d)= 2 \e (s) \left({\frac {qj^*_s}{s}}\right)\ \sqrt{s}\
\left(_q\l_{s}^{d^2}\right)^j $$

\noindent where now we let  $_q\l_{s}$ be the primitive $s$th root of unity
$\xi_s^{-q^*_s \ (\frac{s+1}{2})^2}.$ Thus the rows of $A$ are nonzero
multiples
of  the rows  of a Vandermonde matrix with a non-zero determinant. Here it is
important that as $d$
ranges over $1 \le d \le \frac{s-1}{2},$ $d^2$ does not repeat itself modulo
$s.$

Now we study the matrix  $B.$   For $1 \le d \le s-2,$ and $1 \le j \le s-2,$
with $d$ and $j$ odd , by \cite{LL} we have: $$ G_+(2s,q,\hat d  , j^*_p)=\Cal
G_{2s}(q j_p^*,d) = 2\e (s) \left({\frac {2q j^*_p}{s}}\right)\ \sqrt{s}\
\left(_q\o_{s}^{d^2}\right)^j $$

\noindent where now we let  $_q\o_{s}$ be the primitive $s$th root of unity
$\xi_s^{-q^*_p (\frac{s+1}{2})^3}.$  Let  $_q\theta_{s}$ be the
primitive $s$th root of
unity $_q\o_{s}^2.$
So
$\left(_q\o_{s}^{d^2}\right)^j =
{_q\o_{s}}^{d^2} \left(_q\theta_{s}^{d^2}\right)^\frac{j-1}{2}.$  We may take a
nonzero factor of $_q\o_{s}^{d^2} $ out of each column .
The rows of the resulting matrix are non-zero  multiples of the rows  of a
Vandermonde matrix with a non-zero determinant.
Here it is important that as $d$ ranges over odd numbers such
that $1 \le d \le s-2,$  $d^2$ does not repeat itself modulo $s.$

This completes the proof.

\head \S5 Proof of Theorem 3 \endhead

There are only two lens spaces of order nine up to (not necessarily orientation
preserving) diffeomorphism: $L(9,1),$ and $L(9,4).$ Thus we need only concern
ourselves with $q=1$ and $q=4.$

Using Mathematica we have found that the subset of $\Bbb S(9,1)$ with trivial
WRT invariants is spanned by $(-z^{15} + z^{27})\m_1+(z^{12}-z^{24})\m_2-
z^{15} \m_3+\m_4.$ Similarly the subset of $\Bbb S(9,4)$ with trivial
WRT invariants is
spanned by $(-z^{84}+z^{108}) \m_0 +( z^{60}-z^{72})\m_2 -z^{30}\m_3 +\m_4.$ No
non-zero multiple of these elements lies
in the ordinary skein module over $\L.$

 \enddocument \vfill\eject \end